\documentclass[11pt]{amsproc}
\usepackage{amssymb,amsmath,amsthm,anysize}
\usepackage{longtable}
\usepackage{booktabs}
\usepackage{multirow,bigstrut,bigdelim,color}
\usepackage[shortlabels]{enumitem}
\usepackage{multicol}

\marginsize{2cm}{2cm}{2cm}{2cm}

\newtheorem{Thm}{Theorem}[section]

\newtheorem{Rem}[Thm]{Remark}

\newcommand{\GL}{\textrm{GL}}
\newcommand{\SL}{\textnormal{SL}}

\newcommand{\GU}{\textrm{GU}}

\newcommand{\SU}{\textrm{SU}}

\newcommand{\SO}{\textrm{SO}}
\newcommand{\GO}{\textrm{GO}}
\newcommand{\GSp}{\textrm{GSp}}
\newcommand{\GSO}{\textrm{GSO}}

\newcommand{\Sp}{\textrm{Sp}}

\renewcommand\le{\leqslant}
\renewcommand\ge{\geqslant}

\setlist[description]{leftmargin=1ex,font=\normalfont\bfseries\space,style=nextline}

\begin{document}

\title[Finding involutions with small support]{Finding involutions with small support}

\author{Alice C. Niemeyer and Tomasz Popiel}

\address{
Lehr- und Forschungsgebiet Algebra\\
RWTH Aachen University\\
52062 Aachen\\
Germany.\newline
Email: \texttt{alice.niemeyer@math.rwth-aachen.de}.
}

\address{
Centre for the Mathematics of Symmetry and Computation\\
School of Mathematics and Statistics\\
The University of Western Australia\\
35 Stirling Highway, Crawley, W.A. 6009, Australia.\newline
Email: \texttt{tomasz.popiel@uwa.edu.au}.
}

\thanks{This research forms part of the Australian Research Council Discovery Project DP140100416. 
The second author thanks RWTH Aachen University for financial support and hospitality during his visit in September 2015.}



\maketitle

\begin{abstract}
\noindent
We show that the proportion of permutations $g$ in $S_n$ or $A_n$ such that $g$ has even order and $g^{|g|/2}$ is an involution with support of cardinality at most $\lceil n^\varepsilon \rceil$ is at least a constant multiple of $\varepsilon$. 
Using this result, we obtain the same conclusion for elements in a classical group of natural dimension $n$ in odd characteristic that have even order and power up to an involution with $(-1)$-eigenspace of dimension at most $\lceil n^\varepsilon \rceil$ for a linear or unitary group, or $2\lceil \lfloor n/2 \rfloor^\varepsilon \rceil$ for a symplectic or orthogonal group.
\end{abstract}

\section{Introduction} \label{sec:intro}

Involutions and their centralisers are important not only in the Classification of Finite Simple Groups, but also in many algorithms for computing with finite groups, as can be seen from the numerous applications of Bray's algorithm~\cite{Bray} for finding the centraliser of an involution. 
In practice, involutions (elements of order $2$) can be obtained by raising an element $g$ of even order to the power of $|g|/2$, and this is especially useful in groups where involutions are rare. 
For example, in a finite symmetric group $S_n$ or alternating group $A_n$, algorithms that obtain involutions in this way include those in \cite{BLNPS,BealsPermMod,BP,JLNP}. 
The algorithms in \cite{BLNPS,BealsPermMod,BP} require transpositions, while the algorithm in \cite{JLNP} requires involutions with support of cardinality at most a constant times $\sqrt{n}$ (the {\em support} of a permutation being the set of points that it moves). 
Although involutions in $S_n$ or $A_n$ are rare, it was proved in \cite[Theorem~1.2]{JLNP} that the proportion of elements $g$ in $S_n$ or $A_n$ that have even order and for which the support of $g^{|g|/2}$ has cardinality at most $4\sqrt{n}/3$ is at least $(13 \log(n))^{-1}$. 
Here we strengthen this result considerably. 

\begin{Thm} \label{mainThm}
Let $\varepsilon \in (0,1)$, and let $n$ be an integer satisfying $\lceil (\log(n)+1)^2 \rceil < \lceil n^\varepsilon \rceil \le n - 2\lceil\log(n)\rceil$. 
Let $p(n,\varepsilon)$ denote the proportion of elements $g\in S_n$ such that $g$ has even order and $g^{|g|/2}$ is an involution whose support has cardinality at most $\lceil n^\varepsilon \rceil$, and let $\tilde{p}(n,\varepsilon)$ denote the corresponding proportion of elements in $A_n$. 
Then $p(n,\varepsilon) > \varepsilon/48$ and $\tilde{p}(n,\varepsilon) > \varepsilon/96$.
\end{Thm}

Involutions play an important role in several algorithms for computing with matrix groups or black-box groups; see, for example, the survey article by O'Brien~\cite{Eamonn}.
As an application of Theorem~\ref{mainThm}, we obtain a lower bound for the proportion of elements in a classical group in odd characteristic that have even order and power up to an involution with a $(-1)$-eigenspace of `small' dimension. 
Here, for a finite group $H$ and a subset $I$ of involutions in $H$, we write $P(H, I) = \{ h \in H : |h|\text{ is even and } h^{|h|/2} \in I \}$. 

\begin{table}[t]
\begin{center}
\begin{tabular}{lllll}
\toprule
$n$ & $S$ & $X$ & $\alpha$ & $c_1$ \\ 
\midrule
$\ell$ & $\SL_n(q)$ & $\GL_n(q)$ & $1$ & $1/2$ \\
$\ell$ & $\SU_n(q)$ & $\GU_n(q)$ & $1$ & $1/2$ \\
$2\ell$ & $\Sp_n(q)$ & $\GSp_n(q)$ & $2$ & $1/4$ \\
$2\ell+1$ & $\SO_n(q)$ &  $\GSO_n(q)$ & $2$ & $1/4$ \\
$2\ell$ & $\SO^\pm_n(q)$ & $\GO^\pm_n(q)^\circ$ & $2$ & $1/4$ \\
\bottomrule
\end{tabular}
\end{center}
\caption{Definitions for Theorem~\ref{matrixThm}.}\label{tab:G} 
\end{table}

\begin{Thm} \label{matrixThm}
Let $\varepsilon \in (0,1)$, and let $\ell$ be an integer satisfying $\lceil (\log(\ell)+1)^2 \rceil < \lceil \ell^\varepsilon \rceil \le \ell - 2\lceil\log(\ell)\rceil$. 
Let $q$ be a power of an odd prime, and let $n=n(\ell)$, $S=S(n,q)$, $X=X(n,q)$, $\alpha$ and $c_1$ be as in one of the lines of Table~$\ref{tab:G}$. 
Let $H$ be a group satisfying $S\le  H \le X$, and define $c_2=1/4$ if $S < H \le X$ in lines $3$--$5$ of Table~$\ref{tab:G}$, and $c_2=1$ otherwise.
Define $I$ to be the set of involutions in $H$ that have $(-1)$-eigenspace of dimension $r$ such that $r \le \alpha \lceil \ell^\varepsilon\rceil$. 
Then $|P(\overline{H},\bar{I})|/|\overline{H}| \ge |P(H,I)|/|H| > c_1c_2\cdot\varepsilon/48$, where $\overline{H} = H/Z(H)$ and $\bar{I} = IZ(H)/Z(H)$, with $Z(H)$ the centre of $H$. 
\end{Thm}

Theorems~\ref{mainThm} and~\ref{matrixThm} are proved in Sections~\ref{sec2} and~\ref{sec3}, respectively.

\section{Proof of Theorem~\ref{mainThm}} \label{sec2}

We first prove the result for $S_n$. 
Let $a$ be a positive integer, and $k$ an odd positive integer with $2^ak \le n$.
Consider a permutation $g \in S_n$ with a cycle of length $2^ak$, and the remaining $n-2^ak$ points lying in cycles of lengths not divisible by $2^a$. 
Then the permutation $g^{|g|/2}$ is an involution whose support has cardinality $2^ak$. 
In particular, consider those such permutations for which $a \le A := \log_2(\lceil \log(n) \rceil)$ and $k \le K := \left\lfloor \lceil n^\varepsilon \rceil / \lceil \log(n) \rceil \right\rfloor$, so that $2^ak \le \lceil n^\epsilon \rceil$. 
Noting that the proportion of $(2^ak)$-cycles in $S_n$ is $1/(2^ak)$, and letting $\mathsf{s}_{\neg 2^a}(\ell)$ denote the proportion of elements in $S_{\ell}$ (for some $\ell$) with no cycles of lengths divisible by $2^a$, we have
\begin{equation} \label{Sn-1}
p(n,\varepsilon) \ge \sum_{\substack{k=1 \\ k \text{ odd}}}^K \sum_{a=1}^{A} \frac{\mathsf{s}_{\neg 2^a}(n-2^ak)}{2^ak}.
\end{equation}
Now, $n \ge \lceil n^\varepsilon \rceil + 2\lceil\log(n)\rceil \ge (k+2)\lceil\log(n)\rceil \ge 2^a(k+2)$, so $2^a \le (n-2^ak)/2$. 
Hence, we may apply \cite[Lemma~4.2(a)]{LNP} (a consequence of \cite[Theorem~2.3(b)]{BealsProp}), which gives $\mathsf{s}_{\neg 2^a}(n-2^ak) \ge (4(n-2^ak))^{-1/(2^a)}$. 
Therefore,
\begin{align*}
p(n,\varepsilon) &\ge \frac{1}{4} \sum_{\substack{k=1 \\ k \text{ odd}}}^K \sum_{a=1}^{A} \frac{1}{2^ak(n-2^ak)^{1/2^a}} \\
&> \frac{1}{4} \sum_{\substack{k=1 \\ k \text{ odd}}}^K \sum_{a=1}^{A} \frac{1}{2^ak n^{1/2^a}} 
= \frac{1}{4} \left( \sum_{\substack{k=1 \\ k \text{ odd}}}^K \frac{1}{k} \right) \left( \sum_{a=1}^{A} \frac{1}{2^a n^{1/2^a}} \right).
\end{align*}
Writing $2m+1=k$ and $K' = \lfloor (K-1)/2 \rfloor$, we have
\[
p(n,\varepsilon) > \frac{1}{4} \left( \sum_{m=0}^{K'} \frac{1}{2m+1} \right) 
\left( \sum_{a=1}^{A} \frac{1}{2^an^{1/2^a}} \right).
\]
Since $1/(2m+1)$ is decreasing in $m$ and $(2^a n^{1/2^a})^{-1}$ is increasing in $a$, we can bound the above sums by integrals as follows:
\begin{align}
p(n,\varepsilon) &> \frac{1}{4} \left( \int_0^{K/2} \frac{\text{d}x}{2x+1} \right) \left( \int_0^A \frac{\text{d}x}{2^xn^{1/2^x}} \right) \label{Sn-2} \\
&= \frac{1}{4} \left[\frac{\log(2x+1)}{2}\right]_{0}^{K/2} \left[\frac{1}{\log(2)\log(n)n^{1/2^x}}\right]_{0}^{A} \nonumber \\
&= \frac{\log(K+1)}{8\log(2)\log(n)} \left( \frac{1}{n^{1/2^A}} - \frac{1}{n} \right). \nonumber
\end{align}
Since $K = \left\lfloor \lceil n^\varepsilon \rceil / \lceil \log(n) \rceil \right\rfloor > \lceil n^\varepsilon \rceil / \lceil \log(n) \rceil - 1 \ge n^\varepsilon/(\log(n)+1) - 1$, we have $K+1 > n^\varepsilon/(\log(n)+1)$. 
Therefore,
\[
p(n,\varepsilon) > \frac{1}{8\log(2)} \left( \varepsilon - \frac{\log(\log(n)+1)}{\log(n)} \right) \left( \frac{1}{e} - \frac{1}{n} \right).
\]
Since $\lceil n^\varepsilon \rceil > \lceil (\log(n)+1)^2 \rceil$, we have $n^{\varepsilon/2} > \log(n)+1$, and so the term in the first set of parentheses above is at least $\varepsilon/2$. 
Hence,
\begin{equation} \label{Sn-3}
p(n,\varepsilon) > \frac{\varepsilon}{16\log(2)} \left( \frac{1}{e} - \frac{1}{n} \right).
\end{equation}
Since $n \ge \lceil n^\varepsilon \rceil + 2\lceil \log(n) \rceil > \lceil (\log(n)+1)^2 \rceil + 2\lceil \log(n) \rceil$, we have, in particular, $n\ge 27$. 
Therefore, $1/e - 1/n > \log(2)/3$, and hence $p(n,\varepsilon) > \varepsilon/48$.

The proof for $A_n$ requires only minor changes. 
We again consider permutations with a cycle of length $2^ak$ and the remaining $n-2^ak$ points lying in cycles of length not divisible by $2^a$, but now the product of the cycles not divisible by $2^a$ should lie in $S_n\setminus A_n$. 
Let $\mathsf{a}_{\neg 2^a}(\ell)$ and $\mathsf{c}_{\neg 2^a}(\ell)$ denote, respectively, the proportions of elements in $A_\ell$ and $S_\ell\setminus A_\ell$ (for some $\ell$) with no cycles of lengths divisible by $2^a$. 
Since $\mathsf{c}_{\neg 2^a}(\ell) = 2\mathsf{s}_{\neg 2^a}(\ell) - \mathsf{a}_{\neg 2^a}(\ell)$, \cite[Theorem~3.3(b)]{BealsProp} gives $\mathsf{c}_{\neg 2^a}(\ell) \ge (1-1/(2^a-1))\mathsf{s}_{\neg 2^a}(\ell)$. 
For $a\ge 2$, we therefore have $\mathsf{c}_{\neg 2^a}(\ell) \ge 2/3 \cdot \mathsf{s}_{\neg 2^a}(\ell)$, and so the inequality \eqref{Sn-1} is replaced by
\[
\tilde{p}(n,\varepsilon) \ge \frac{2}{3} \sum_{\substack{k=1 \\ k \text{ odd}}}^K \sum_{a=2}^{A} \frac{\mathsf{s}_{\neg 2^a}(n-2^ak)}{2^ak}.
\]
The range of the second integral in \eqref{Sn-2} is accordingly changed to $[1,A]$, and so instead of \eqref{Sn-3} we obtain
\[
\tilde{p}(n,\varepsilon) > \frac{2}{3} \cdot \frac{\varepsilon}{16\log(2)} \left( \frac{1}{e} - \frac{1}{\sqrt{n}} \right).
\]
Since $n \ge 27$, we have $1/e-1/\sqrt{n} > \log(2)/4$, and hence $\tilde{p}(n,\varepsilon) > \varepsilon/96$.

\begin{Rem}
\textnormal{
If the single cycle of length $2^ak$ in the proof of Theorem~\ref{mainThm} is replaced by a product of cycles of lengths divisible by $2^a$ but not by $2^{a+1}$, then the resulting permutation $g$ still powers up to an involution with support at most $\lceil n^\varepsilon \rceil$. 
By \cite[Lemma~5.2]{JLNP} (a refinement of a particular case of \cite[Theorem~1.1]{Sukru}), the proportion $t_{2^a}(2^ak)$ of elements in $S_{2^ak}$ with all cycles of lengths divisible by $2^a$ but not by $2^{a+1}$ is at least $c(a)/k^{1-1/2^{a+1}}$, where $c(a) = 1/(2^{2a}3^{1/2^{a+1}})$. 
For a fixed value of $a$, this lower bound, viewed as a function of $k$, is asymptotically better than the lower bound of $1/(2^ak)$ that we have employed in the proof of Theorem~\ref{mainThm}, and it was used in the proof of \cite[Theorem~1.2]{JLNP}, where only two values of $a$ were considered. 
However, in the proof of Theorem~\ref{mainThm}, where we sum over a range of values of $a$ that depends on $n$, the lower bound $c(a)/k^{1-1/2^{a+1}}$ for $t_{2^a}(2^ak)$ is too small to yield a constant lower bound for $p(n,\varepsilon)$, and so we have opted to use the `weaker' bound $1/(2^ak)$.
}
\end{Rem}

\section{Proof of Theorem~\ref{matrixThm}} \label{sec3}

The proof is essentially the same as that of \cite[Theorem~1.1]{LNP}, but we change the range of $r$ from $r \in [n/3,2n/3)$ to $r \le \alpha \lceil \ell^\varepsilon \rceil$. 
Consider first the case where $H = \GL_n(q)$ or $\GU_n(q)$, and let $W$ be the Weyl group of $H$, noting that $W \cong S_n$. 
We modify the definition of $M(a)\subset W$ from \cite[Section~5(ii)]{LNP} as follows: an element $w \in W$ lies in $M(a)$ if and only if it contains a cycle of length $2^a k$ for some odd integer $k$ with $2^a k \le \lceil \ell^\varepsilon\rceil = \lceil n^\varepsilon\rceil$, and no other cycle has length divisible by $2^a$. 
The argument in \cite[Section~5(ii)]{LNP} is then still valid, and hence we still obtain \cite[Equation~(5)]{LNP}, which gives $|P(H,I)|/|H| \ge p(n,\varepsilon)/2$, showing that $c_1 = 1/2$. 
If $H=\Sp_{2\ell}(q)$, $\SO_{2\ell+1}(q)$ or $\SO^\pm_{2\ell}(q)$, then the argument in \cite[Section~(vii)]{LNP} still applies: the Weyl group is now a subgroup of $S_2 \wr S_\ell$, $M(a)$ is redefined so that the cycle of length $2^ak \le \lceil \ell^\varepsilon \rceil$ is positive, and we conclude that $c_1 = 1/4$. 
The same arguments as in \cite[Sections~(v) and~(viii)]{LNP} give the claimed values for $c_2$. 
Finally, the proof of \cite[Corollary~1.2]{LNP} gives the result for the groups $\overline{H}$.

\end{document}